\documentclass[a4paper,12pt]{article}
\usepackage{graphicx,amssymb,amstext,amsmath}

\usepackage[latin1]{inputenc}
\usepackage{amsthm}
\usepackage{dsfont}
\usepackage{amsmath}
\usepackage{amsfonts}
\usepackage{amssymb}
\usepackage{mathrsfs}
\usepackage{graphicx}
\usepackage{cancel}
\usepackage{comment}
\usepackage{calligra}
\usepackage{graphicx,color}
\usepackage[all]{xy}
\usepackage[total={5.2in,9.1in},
top=1.4in, left=1.55in, includefoot]{geometry}

\usepackage{amsmath}
\usepackage{amssymb}
\usepackage{amsthm}
\newtheorem{teo}{Theorem}[section]
\newtheorem{lem}[teo]{Lemma}
\newtheorem{prop}[teo]{Proposition}
\newtheorem{coro}[teo]{Corollary}
\newtheorem{defi}[teo]{Definition}

\newtheorem{remark}[teo]{Remark}

\newtheorem{conj}{Conjecture}

\begin{document}
\thispagestyle{empty}

\title{\textbf{\Large{Curvature Rigidity Through Level Sets of Lyapunov Exponents in Geodesic Flows}}}
\author{ Sergio Roma\~na}
\date{}
\maketitle

\begin{abstract}
In this paper, we establish new geometric rigidity results through the study of Lyapunov exponent level sets via invariant measures. First, we 
prove that for a manifold $M$ without focal points, if the zero Lyapunov exponent level set has full measure with respect to a fully supported 
invariant measure, then $M$ must be flat. This result recovers and extends a result of Freire and Ma\~n\'e (cf. \cite{MF}). 

Second, we prove that if the level set of the Lyapunov exponents has full measure with respect to some fully supported measure, then the sectional curvature must be constant. This advances the resolution of \cite[Conjecture 1]{IR2}.

Furthermore, we establish curvature relationships between manifolds with $1$-equivalent geodesic flows, yielding a new criterion to obstruct smooth conjugacy for flows on manifolds without conjugate points.
Our techniques provide unified proofs for all rigidity results in \cite{IR2} as corollaries, and additionally yield rigidity theorems for totally geodesic submanifolds in settings without conjugate points. Notably, several key results hold without compactness assumptions.

\begin{quote}
\end{quote}

\end{abstract}

\section{Introduction}
The rigidity of Lyapunov exponents has become a highly interesting topic for exploring geometric properties of manifolds by studying the asymptotic behavior of the derivative of geodesic flows. Many authors have addressed these rigidity problems in the context of geodesic flows, as well as for different types of dynamics. 
For a comprehensive overview of the state of the art, the reader may refer to the introductions of \cite{Bu, IR2, Nina2}.\\
In this paper, we study the rigidity of Lyapunov exponents from the perspective of invariant measures for the geodesic flow on manifolds without conjugate points. This approach enables progress toward resolving the following conjecture, proposed in \cite{IR2} and inspired by \cite{Bu}:

\begin{conj}\label{conj.1}
Let $M$ be a complete Riemannian manifold of finite volume whose geodesic flow is Anosov. If the unstable Lyapunov exponents are constant across all periodic orbits, then $M$ has constant negative sectional curvature.
\end{conj}
This conjecture was proved in dimension two by the author jointly with N. Nina (cf. \cite{Nina2}), while C. Butler established the compact case of negative curvature (cf. \cite{Bu}).

We further explore specific invariant measures, such as the Liouville measure, measures of full support, and periodic measures, to establish curvature rigidity. In particular, we analyze the measure structure of the level sets of non-negative Lyapunov exponents:
$$\Lambda_\alpha=\{\theta\in SM: \chi^{+}(\theta)=\alpha(n-1)\},$$
where $\alpha\geq 0$.\\
Our first result recovered a result by Freire and Ma\~n\'e (see. \cite{MF}) and introduces a new equivalence involving the measure of  $\Lambda_0$.

\begin{teo}\label{Theorem_Flat}
    Let $M$ be a compact manifold without focal points. Then the following are equivalent:
    \begin{itemize}
        \item[\emph{(a)}] $h_{\text{top}}(\phi^t_{M})=0$.
        \item[\emph{(b)}] There exists a fully supported invariant measure $\mu$ with $\mu(\Lambda_0)=1$.
        \item[\emph{(c)}] $M$ is flat.
    \end{itemize}
\end{teo}   

In dimension two, Theorem \ref{Theorem_Flat} remains valid under the weaker assumption of no conjugate points (see \cite{hopf}). However, the extension of this result to higher dimensions under only the no conjugate points assumption remains an open problem (cf. \cite{MF}).

For surfaces of non-positive curvature, we have
\[
\Lambda_0 = \{\theta \in SM : \gamma_\theta(t) = 0 \text{ for all } t \in \mathbb{R}\},
\]
commonly known as the \emph{irregular set}. A central conjecture states that this set has zero Liouville measure, which would imply ergodicity of the geodesic flow for such surfaces (see \cite{KB}). While the ergodicity question remains unresolved, it is established that the geodesic flow $\phi^t_M$ restricted to $\Lambda_0$ has zero topological entropy (cf. \cite{WU}).

As a consequence of Theorem \ref{Theorem_Flat}, we obtain the following more general result in higher dimensions.

\begin{coro}\label{Zero_entropy}
For any compact manifold $M$ without focal points, the topological entropy satisfies $h_{\mathrm{top}}(\phi^t_M, \Lambda_0) = 0$.
\end{coro}

Our next result explores rigidity for positive values of $\alpha$ from the perspective of invariant measures supported on $\Lambda_\alpha$, requiring only the no conjugate points condition.

\begin{teo}\label{Theorem_Rigidity}
Let $M$ be a manifold without conjugate points and with curvature bounded below. If $\mu$ is an invariant measure satisfying $\mu(\Lambda_{\alpha}) = 1$ for some $\alpha \geq \sqrt{-\int_{SM} \text{Ric}\,d\mu}$, then the sectional curvature $K$ is constant on $\pi(\text{Supp}(\mu))$ and equals $\int_{SM} \text{Ric}(\theta)\,d\mu$.
\end{teo}

When $\mu$ in Theorem \ref{Theorem_Rigidity} is ergodic and has full support, the curvature must be constant everywhere. More precisely:

\begin{coro}\label{C_Theorem_Rigidity}
If $\mu$ is an ergodic measure with full support satisfying $\mu(\Lambda_{\alpha}) > 0$ for some $\alpha \geq \sqrt{-\int_{SM} \text{Ric}(\theta)\,d\mu}$, then the sectional curvature $K$ is constant and equals $\int_{SM} \text{Ric}(\theta)\,d\mu$.
\end{coro}
The Theorem \ref{Theorem_Rigidity} and Corollary \ref{C_Theorem_Rigidity} extended the Rigidity on \cite[Theorem 1.1]{IR2} (see Section \ref{Some Applications}).

Our next result establishes a proof of Conjecture \ref{conj.1} for values of Lyapunov exponents in periodic orbits that are greater than or equal to $\sqrt{-\int_{SM} \text{Ric}(\theta)\,d\mathcal{L}}$. Specifically:

\begin{teo}\label{Theorem_Rigidity2}
Let $M$ be a compact manifold with Anosov geodesic flow $\phi^t_M$. If the positive Lyapunov exponent is constant on periodic orbits with value $\alpha$, where $\alpha \geq \sqrt{-\int_{SM} \text{Ric}(\theta)\,d\mathcal{L}}$, then $M$ has constant sectional curvature.
\end{teo}

From the perspective of periodic measures, we can obtain rigidity results even in the non-compact case:

\begin{teo}\label{Theorem_Rigidity3}
Let $M$ be a complete manifold without conjugate points such that $\phi^t|_{\overline{\text{Per}(\phi^t)}}$ is ergodic. If
\[
\chi^+(\theta) \geq (n-1)\sqrt{-\int_{SM} \text{Ric}(\theta)\,d\mu_{\theta}} \quad \text{for all } \theta \in \text{Per}(\phi^t),
\]
then $M$ has constant sectional curvature on $\pi(\overline{\text{Per}(\phi^t)})$.
\end{teo}

For Anosov geodesic flows on finite volume manifolds, we know that periodic orbits form a dense set (see \cite{Nina}) and the flow is ergodic. This allows us to prove Conjecture \ref{conj.1} for certain values of Lyapunov exponents. More precisely:

\begin{teo}\label{C_Theorem_Rigidity3}
Let $M$ be a complete finite-volume manifold with Anosov geodesic flow. If the positive Lyapunov exponents on periodic orbits are constant and equal to $\alpha$, where
\[
\alpha \geq \sqrt{-\int_{SM} \text{Ric}(\theta)\,d\mu_{\theta}} \quad \text{for all } \theta \in \text{Per}(\phi^t),
\]
then $M$ has constant sectional curvature.
\end{teo}

The second part of this paper deals with regular equivalences between geodesic flows of manifolds without conjugate points. 

We say that two geodesic flows $\phi^t_{M}$ and $\phi^t_{N}$ are \emph{equivalent} if there exists a continuous map $h\colon SM\to SN$ and a time reparametrization $s\colon \mathbb{R} \to \mathbb{R}$ such that 
\[ h \circ \phi^t_{M} = \phi_{N}^{s(t)} \circ h \quad \text{for all } t\in \mathbb{R}. \]
Here, $h$ is called an \emph{equivalence} and $s$ a \emph{reparametrization}. 

When $h$ is a homeomorphism, we call it a \emph{conjugacy}. An equivalence (or conjugacy) $h$ is called a \emph{$1$-equivalence} (respectively \emph{$1$-conjugacy}) if $h$ is bi-Lipschitz, i.e., there exist constants $C_1, C_2 > 0$ such that
\begin{equation}\label{eq:bi_lipschitz}
C_1\cdot d(x,y) \leq D(h(x),h(y)) \leq C_2\cdot d(x,y),
\end{equation}
where $d$ and $D$ denote the distances on $SM$ and $SN$, respectively.

If $s(t) \geq t$, $t\geq 0$ and $h$ is a $1$-equivalence, we say that $\phi^t_{M}$ and $\phi^t_{N}$ satisfy property $\mathcal{E}_1$. For a comprehensive discussion of flow equivalences with $s(t) = t$, we refer to \cite{IR2}.

\begin{remark}
Since $1$-equivalences are bi-Lipschitz, we have $\dim M \leq \dim N$. Moreover, there exists a set $\Gamma_{M,N}\subset SM$ of full Liouville measure where $h$ is differentiable. Our next result establishes curvature rigidity using invariant measures supported on $\Gamma_{M,N}$.
\end{remark}

\begin{teo}\label{N-Conjugacy} 
Let $M$ and $N$ be complete Riemannian manifolds whose geodesic flows $\phi_{M}^{t}$ and $\phi_{N}^{t}$ satisfy property $\mathcal{E}_1$ with $1$-equivalence $h$. Assume that: $M$ and $N$ have curvature bounded below, $M$ has no conjugate points, and there exists an invariant measure $\mu$ with $\mu(SM\setminus \Gamma_{M, N}) = 0$.
Then 
\[ \sup K_{N} \geq \int_{SM} \text{Ric}(\theta) d\mu. \]
Equality holds if and only if one of the following conditions is satisfied:
\begin{itemize}
    \item[\emph{\text{1.}}] $K_{M}|_{\pi(\text{Supp}(\mu))} = \sup K_N$; 
    \item[\emph{\text{2.}}] $\displaystyle\lim_{t\to +\infty}\frac{s(t)}{t} = 1$.
\end{itemize}
Furthermore, if $\dim M = \dim N$, then $N$ has constant sectional curvature on $\pi(\text{Supp}(h_{*}(\mu)))$.
\end{teo}

We know that the \emph{Liouville measure} $\mathcal{L}$ satisfies $\mathcal{L}(SM\setminus \Gamma_{M,N})=0$. Therefore, under the same hypotheses as Theorem \ref{N-Conjugacy}, we obtain the following corollary:

\begin{coro}\label{N-Conjugacy-Liouville}
Under the same hypotheses as \emph{Theorem \ref{N-Conjugacy}}, the following statements are equivalent:
\begin{itemize}
\item $\displaystyle\sup K_{N} = \int_{SM} \text{Ric}(\theta)\,d\mathcal{L}$
    \item $K_M \equiv \sup K_N$;
    \item $\displaystyle\lim_{t\to +\infty}\frac{s(t)}{t}=1$.
\end{itemize}
Moreover, if $\dim M = \dim N$, then $N$ has constant sectional curvature.
\end{coro}

The Theorem \ref{N-Conjugacy} and Corollary \ref{N-Conjugacy-Liouville} extended \cite[Theorem 1.2 and Theorem 1.4]{IR2} (see Section \ref{Some Applications}).

When the $1$-equivalence is actually $C^1$, we have $\Gamma_{M,N}=SM$, allowing us to characterize rigidity through periodic points. Specifically, if $\mu_\theta$ denotes the invariant measure supported on a periodic orbit $\theta$, then:

\begin{teo}\label{C1-Conjugacy}
Let $M$ and $N$ be complete Riemannian manifolds whose geodesic flows $\phi_{M}^{t}$ and $\phi_{N}^{t}$ satisfy property $\mathcal{E}_1$ with $C^1$-equivalence $h$. If $M$ has finite volume and $\phi^t_{M}$ is Anosov, then
\[
\sup\,K_{N}|_{\mathcal{O}(h(\theta))} \geq \int_{SM} \text{Ric}(\theta)\,d\mu_{\theta} \quad \text{for all } \theta \in \text{Per}(\phi^t_{M}).
\]
Equality holds for all geodesics if and only if:
\begin{itemize}
    \item $M$ has constant sectional curvature;
    \item $\displaystyle\lim_{t\to +\infty}\frac{s(t)}{t}=1$.
\end{itemize}
Furthermore, if $\dim M = \dim N$, then $N$ has constant sectional curvature.
\end{teo}

\begin{remark}
If the inequality in \emph{Theorem \ref{C1-Conjugacy}} is strict for at least one periodic orbit, then the geodesic flows of $N$ and $M$ cannot be $C^1$-equivalent. This provides an effective criterion for determining when two geodesic flows fail to be $C^1$-equivalent.
\end{remark}

\section{Preliminaries and Notation}
Throughout this paper, $M$ and $N$ will denote complete, boundaryless Riemannian manifolds without boundary with $2\leq \text{dim} \,M\leq \text{dim}\, N$. We will use $TM$ and $TN$ to represent their tangent bundles, $SM$ and $SN$ to denote their unit tangent bundles. The map $\pi\colon SM \to M$ will denote the canonical projection and  $\mathcal{L}$ the \emph{Liouville} measure on $SM$.\\

We offer a brief introduction to geodesic flows and their geometric properties. For a more in-depth discussion on topics such as geodesic flows, Anosov geodesic flows, no conjugate points, no focal points, Green sub-bundles, and other related concepts, we refer the reader to \cite{IR1, IR2, Nina, CR, Kn}.

\subsection{Geodesic Flow}
For each $t \in \mathbb{R}$ and $\theta = (x, v) \in SM$, the family of diffeomorphisms $\phi^t_{M}\colon SM \to SM$ is defined by
$\phi^t_{M}$
\[
\phi^t_{M}(\theta) = (\gamma_{\theta}(t), \gamma_{\theta}'(t)),
\]
where $\gamma_{\theta}(t)$ is the unique geodesic with initial conditions $\gamma_{\theta}(0) = p$ and $\gamma_{\theta}'(0) = v$. This family is called the \emph{geodesic flow} on $SM$.

The spray vector field $G$ is the derivative of $\phi^t$, given by
\[
G(\theta) := \frac{d}{dt} \phi^t_{M}(\theta)\Big|_{t=0}.
\]

For $\theta = (x, v) \in TM$, we identify $T_{\theta}TM$ with $T_xM \oplus T_xM$ via the map $\xi \to (d\pi_{\theta}(\xi), \mathcal{K}_{\theta}(\xi))$, where $\pi\colon TM \to M$ is the canonical projection, and $\mathcal{K}\colon TTM \to TM$ is the connection map defined by the Levi-Civita connection. This identification allows us to define the \emph{Sasaki metric} on $TM$ as:
\[
\langle \xi, \eta \rangle_{\theta} = \langle d \pi_{\theta} (\xi), d \pi_{\theta} (\eta) \rangle + \langle \mathcal{K}_{\theta}(\xi), \mathcal{K}_{\theta}(\eta) \rangle.
\]

The geodesic flow $\phi^t:SM \rightarrow SM$ is Anosov (with {respect to} the Sasaki metric on $SM$) if $T(SM)$ has a splitting $T(SM) = E^s \oplus \langle G \rangle \oplus E^u $ such that 
\begin{eqnarray*}
	d\phi^t_{\theta} (E^s(\theta)) &=& E^s(\phi^t(\theta)),\\
	d\phi^t_{\theta} (E^u(\theta)) &=& E^u(\phi^t(\theta)),\\
	||d\phi^t_{\theta}\big{|}_{E^s}|| &\leq& C \lambda^{t},\\
	||d\phi^{-t}_{\theta}\big{|}_{E^u}|| &\leq& C \lambda^{t},
	\end{eqnarray*}
for all $t\geq 0$ with $ C > 0$ and $0 < \lambda <1$,  where $G$ is the vector field derivative of the geodesic flow. The constant $\lambda$ is called  \emph{a constant of contraction}.

Natural examples of Anosov geodesic flows occur on compact manifolds with negative curvature (cf. \cite{Anosov}) and non-compact manifolds with pinched negative curvature (\cite{Kn}). Other examples exist without curvature restrictions (see \cite{Ebe} and \cite{IR1}). For geometric conditions guaranteeing the Anosov property in both compact and non-compact settings, we refer the reader to \cite{IR1, CR2}.

\subsection{Geometric and Dynamical Properties of Conjugate and Focal Points}\label{sec-ncp} 
The dynamics of geodesic flows are intimately connected to the behavior of Jacobi fields. For any geodesic $\gamma$ in a Riemannian manifold $M$, a Jacobi field $J$ along geodesics $\gamma$ satisfies $J'' + R(\gamma', J)\gamma' = 0$, using the horizontal and vertical decomposition. Jacobi fields encode the dynamics of the geodesic flow via $d\phi^t_\theta(\xi) = (J_\xi(t), J'_\xi(t))$, where $J_{\xi}(t)$ denotes the unique Jacobi field along $\gamma_{\theta}(t)$ such that $J_{\xi}(0) = d\pi_{\theta}(\xi)$ and $J'_{\xi}(0) = \mathcal{K}_{\theta}(\xi)$ (see \cite{Pa} for more details).
\noindent A manifold has \emph{no conjugate points} if no geodesic has non-trivial Jacobi fields vanishing at two points, and \emph{no focal points} if $\|J\|^2$ grows strictly for $t > 0$ when $J(0) = 0$, $J'(0) \neq 0$. The latter implies the former.  

For Anosov flows on compact/finite-volume manifolds, the absence of conjugate points was proven in \cite{Klin, Ma}, later extended to infinite volume cases under bounded curvature assumptions in \cite{IR3}. The \emph{Green subbundles} $G^{s/u}_\theta \subset T_\theta SM$, defined for manifolds without conjugate points, consist of vectors whose Jacobi fields remain bounded as $t \to \pm \infty$. If $\dim M = n$, then $\dim G^{s/u}_\theta = n-1$. Continuous dependence on $\theta$ holds when $M$ has no focal points (\cite{Ebe}), and the Anosov condition reduces to $G^s_\theta \cap G^u_\theta = \{0\}$ and $G^*_{\theta} = E^{*}_{\theta}$, where $ * = s, u$ (see \cite{Ebe} for more details on Green subbundles).


\subsubsection{Jacobi Tensors and Riccati Equation}\label{Jacobi Tensors}
Let $\gamma_\theta$ be a geodesic with parallel orthonormal frame $\{V_1(t),...,V_{n-1}(t),\gamma_\theta'(t)\}$. The Jacobi equation for orthogonal vector fields takes the matrix form:
\[
Y''(t) + R(t)Y(t) = 0,
\]
where $R(t)$ is the curvature endomorphism on $N(\phi^t(\theta)) := \gamma_\theta'(t)^\perp$.

The \emph{Riccati equation} for the geodesic flow is:
\begin{equation}\label{Riccati}
U'(t) + U^2(t) + R(t) = 0,
\end{equation}
where $U(t):N(\phi^t(\theta))\to N(\phi^t(\theta))$ is self-adjoint.

The unstable and stable solutions $U_{\theta,u}(t)$ and $U_{\theta,s}(t)$ arise from Jacobi tensors:
\[
U_{\theta,*}(t) = {Y}'_{\theta,*}(t)Y_{\theta,*}^{-1}(t),\quad *\in\{u,s\}
\]
where $Y_{\theta,*}(t)$ are the unstable/stable Jacobi tensors.\\

When the curvature is bounded below by $-c^2$, then (see \cite{Green} and \cite{Ebe})
\begin{equation}\label{EQ-Green}
\max\{\|U_{\theta, s}(t)\|, \|U_{\theta, u}(t)\|\}\leq c.
\end{equation}

When the sectional curvature is negative bounded above, that is, there is a positive constant $b$ such that $K\leq -b^2$, then (cf. \cite{Ebe} and \cite{Kn})
\begin{equation}\label{ec4.8}
\langle U_{\theta,u}(t)x,x\rangle\geq b\langle x,x\rangle
\end{equation}
and
\begin{equation*}
\langle U_{\theta,s}(t)x,x\rangle\leq -b\langle x,x\rangle,
\end{equation*}
for all $t\in \mathbb{R}$.

\subsubsection{Lyapunov Exponents}
The symmetric solutions $U_{\theta,u}(t)$ and  $U_{\theta,s}(t)$ encode all Lyapunov exponents information for the Green subbundles. Precisely, given $\theta\in SM$ and $\xi \in T_{\theta}SM$, the associated \emph{Lyapunov exponent} is defined as 
\[
\chi(\theta, \xi) := \lim_{t \to \infty} \frac{1}{t} \log \|D\phi^t_{\theta}(\xi)\|,
\]
when the limit exists.

For any invariant probability measure $\mu$, Oseledets' Theorem \cite{Osel} ensures the existence of Lyapunov exponents on a full $\mu$-measure set $\Lambda \subset SM$. In the case of Anosov geodesic flows on finite-volume / non-compact manifolds, this existence result extends via \cite{CR}.

The Lyapunov exponents \(\chi(\theta, \cdot)\) are represented as \(\chi^u(\theta, \cdot)\) and \(\chi^s(\theta, \cdot)\) for the unstable and stable Green subbundles, respectively. When \(M\) has no conjugate points, \(\chi^u(\theta, \cdot)\) is non-negative, while \(\chi^s(\theta, \cdot)\) is non-positive.

Moreover, for manifolds without conjugate points (see \cite{MF, IR2})
\begin{equation}\label{EQ-Form-Lyp-Exp}
\chi^{+}(\theta):=\lim_{t \to \infty} \frac{1}{t} \log |\det|D\phi^t_{\theta}|_{G^u_{\theta}}| = \lim_{t \to \infty} \frac{1}{t} \int_{0}^{t} \mathrm{tr}\,U_{\theta,u}(s)ds.
\end{equation}

\subsubsection{Ricci Curvature and Lyapunov Exponents}

The \emph{Ricci Curvature} in a point $(x,v)\in SM$  is defined as the following average (see \cite{doCarmo})
$$\text{Ric}(x,v)=\frac{1}{n-1}\sum_{i=1}^{n-1}R(v,v_i,v, v_i),$$
where $\{v, v_1, \dots, v_{n-1}\}$ is a orthonormal basis of $T_{x}M$.
 
So, using the notation of the Section \ref{Jacobi Tensors}, this becomes: 
$$\text{Ric}(\phi^t_{M}(\theta))=\frac{\text{tr}R(t)}{n-1},$$
where $\text{tr}R(t)$ denotes trace of $R(t)$.\\
From Cauchy-Schwarz inequality, the symmetry of $U_{\theta,u}(t)$,  (\ref{EQ-Green}), and (\ref{EQ-Form-Lyp-Exp}) we can take the trace and integration of the Riccati equation (\ref{Riccati}) yields the fundamental inequality   
$$\chi^{+}(\theta):=\lim_{t \to \infty} \frac{1}{t} \log |\det|D\phi^t_{\theta}|_{G^u_{\theta}}|\leq (n-1)\sqrt{-\lim_{t\to \infty}\frac{1}{t}\int_{0}^{t}\text{Ric}(\phi^s_{M}(\theta))ds}.$$
as established in \cite[Lemma 3.6]{IR2}.

\section{Estimative of Lyapunov Exponents}
For each $\theta\in SM$ we denote the asymptotic upper and lower Ricci averages:  $$\Gamma_{R}^{+}(\theta)= \displaystyle \limsup_{t\to +\infty}\frac{1}{t}\int_{0}^{t}\text{Ric}(\phi^t(\theta))ds \,\,\text{and}\,\, 
\Gamma_{R}^{-}(\theta)= \displaystyle \liminf_{t\to \infty}\frac{1}{t}\int_{0}^{t}\text{Ric}(\phi^t(\theta))ds$$
By the Birkhoff Ergodic Theorem, for any invariant measure $\mu$, these limits coincide $\mu$-almost everywhere, \emph{i.e.},  $$\Gamma_{R}^{+}(\theta)=\Gamma_{R}^{-}(\theta), \,\,\text{for}\,\, \mu\text{-a.e.} \,\, \theta\in SM.$$

\begin{lem}\label{L1NEW}
Let $M$ be a manifold without focal points and curvature bounded below by $-c^2$. For any $\theta\in SM$ where $\chi^{+}(\theta)$ exists,
\begin{equation}\label{EQ0NEW}
-(n-1)\frac{\Gamma_{R}^{+}(\theta)}{c} \leq \chi^{+}(\theta) \leq (n-1)\sqrt{-\Gamma_{R}^{-}(\theta)}.
\end{equation}
\end{lem}

\begin{proof}
Let $U_{\theta,u}(t)$ be the unstable Riccati solution associated to the unstable Jacobi tensor. Since $M$ has no focal points and curvature bounded below by $-c^2$, the symmetric operator $U_{\theta,u}(t)$ has eigenvalues satisfying $0 \leq \lambda_1(t) \leq \cdots \leq \lambda_{n-1}(t) \leq c$. This implies:
\begin{equation}\label{EQ1NEW}
\mathrm{tr}\,U_{\theta,u}^2(t) = \sum_{i=1}^{n-1} \lambda_i^2(t) \leq c\sum_{i=1}^{n-1} \lambda_i(t) = c\,\mathrm{tr}\,U_{\theta,u}(t).
\end{equation}

Taking the trace of the Riccati equation \eqref{Riccati} and integrating yields:
\begin{align*}
0 &= \frac{1}{t}\int_{0}^{t} \mathrm{tr}\,U'_{\theta,u}(s)\,ds + \frac{1}{t}\int_{0}^{t} \mathrm{tr}\,U^2_{\theta,u}(s)\,ds + \frac{n-1}{t}\int_{0}^{t} \mathrm{Ric}(\phi^{s}(\theta))\,ds \\
&= \frac{\mathrm{tr}\,U_{\theta,u}(t) - \mathrm{tr}\,U_{\theta,u}(0)}{t} + \frac{1}{t}\int_{0}^{t} \mathrm{tr}\,U^2_{\theta,u}(s)\,ds + (n-1)\frac{1}{t}\int_0^t \mathrm{Ric}(\phi^s(\theta))\,ds.
\end{align*}

From \eqref{EQ-Green}, \eqref{EQ-Form-Lyp-Exp}, and (\ref{EQ1NEW}) we obtain the lower bound:
\begin{align*}
-(n-1)\Gamma_{R}^{+}(\theta) &= \liminf_{t\to +\infty} \frac{1}{t}\int_{0}^{t} \mathrm{tr}\,U^2_{\theta,u}(s)\,ds \\
&\leq c \liminf_{t\to +\infty} \frac{1}{t}\int_{0}^{t} \mathrm{tr}\,U_{\theta,u}(s)\,ds = c\,\chi^{+}(\theta),
\end{align*}
proving the left inequality in \eqref{EQ0NEW}.

For the upper bound, we apply the Cauchy-Schwarz inequality and the trace inequality for symmetric operators:
\begin{equation*}
(\mathrm{tr}\,U_{\theta,u}(s))^2 \leq (n-1)\,\mathrm{tr}\,U^2_{\theta,u}(s)).
\end{equation*}
This gives:
\begin{align*}
\frac{1}{t}\int_{0}^{t} \mathrm{tr}\,U_{\theta,u}(s)\,ds &\leq \sqrt{\frac{1}{t}\int_{0}^{t} (\mathrm{tr}\,U_{\theta,u}(s))^2\,ds} \\
&\leq \sqrt{\frac{n-1}{t}\int_{0}^{t} \mathrm{tr}\,U^2_{\theta,u}(s))\,ds}.
\end{align*}
Taking limits and using \eqref{EQ-Green}, \eqref{EQ-Form-Lyp-Exp} yields:
\begin{equation}\label{EQLN1}
\chi^{+}(\theta) \leq (n-1)\sqrt{-\Gamma_{R}^{-}(\theta)},
\end{equation}
completing the proof.
\end{proof}

\begin{remark}\label{RNEW1}
The conclusions of \emph{Lemma~\ref{L1NEW}} remain valid for non-compact manifolds. Furthermore, we emphasize that:
\begin{itemize}
    \item The upper bound in inequality~\eqref{EQ0NEW} requires only the absence of conjugate points and and curvature bounded below.
    \item The lower bound additionally requires no focal points and curvature bounded below.
\end{itemize}
This distinction becomes particularly relevant when studying rigidity questions in finite/infinite-volume settings.
\end{remark}

\noindent The Ruelle inequality yields the following lemma, which extends \cite[Corollary II.1]{MF} to the non-compact case:

\begin{lem}\label{RIGI}
Let $M$ be a manifold without conjugate points where the Ruelle inequality holds. Then for any invariant measure $\mu$:
\[ h_{\mu}(\phi^t) \leq (n-1)\sqrt{-\int_{SM}\mathrm{Ric}(\theta)\,d\mu}. \]
Equality occurs if and only if $K|_{\pi(\mathrm{Supp}(\mu))}$ is constant.
\end{lem}

\begin{proof}
Combining the Ruelle inequality with \eqref{EQ0NEW}, we obtain:
\begin{align*}
h_{\mu}(\phi^t) &\leq \int_{SM} \chi^+(\theta)\,d\mu \leq (n-1)\int_{SM} \sqrt{-\widetilde{\mathrm{Ric}}(\theta)}\,d\mu \\
&\leq (n-1)\sqrt{-\int_{SM} \widetilde{\mathrm{Ric}}(\theta)\,d\mu} = (n-1)\sqrt{-\int_{SM} \mathrm{Ric}(\theta)\,d\mu},
\end{align*}
where $\widetilde{\mathrm{Ric}}$ denotes the Birkhoff average of $\mathrm{Ric}$.

\noindent When equality holds, \eqref{EQLN1} implies:
\[ \int_{SM} \mathrm{tr}\,U_u\,d\mu = \int_{SM} \widetilde{\mathrm{tr}\,U_u}\,d\mu = (n-1)\sqrt{-\int_{SM} \mathrm{Ric}(\theta)\,d\mu}. \]
From the Cauchy-Schwarz inequality and Riccati equation  \eqref{Riccati}
\begin{equation}\label{EQ5NEW}
\int_{SM} \mathrm{tr}\,U_u\,d\mu \leq \sqrt{n-1}\sqrt{\int_{SM} \mathrm{tr}\,U_u^2\,d\mu} = (n-1)\sqrt{-\int_{SM} \mathrm{Ric}(\theta)\,d\mu}.
\end{equation}

\noindent Equality in \eqref{EQ5NEW} requires $U_u$ to be a scalar multiple of $I$, which by \eqref{Riccati} forces constant sectional curvature on $\pi(\mathrm{Supp}(\mu))$.
\end{proof}

\begin{remark}\label{R_RIGI}
The rigidity conclusion holds under the weaker assumption:
\[ \int_{SM} \mathrm{tr}\,U_u\,d\mu = (n-1)\sqrt{-\int_{SM} \mathrm{Ric}(\theta)\,d\mu}. \]
In this case, $U_u = \lambda I$ $\mu$-a.e., where $\lambda = \sqrt{-\int_{SM} \mathrm{Ric}(\theta)\,d\mu}$.
\end{remark}

The results of \cite{CR} establish that the Ruelle inequality holds for Anosov geodesic flows when both the curvature tensor and its derivatives are uniformly bounded, making Lemma \ref{RIGI} applicable in this setting.


\begin{lem}\label{L2'NEW}
    Let $M$ be a manifold without focal points and $\mu$ an invariant measure satisfying $\mu(\Lambda_0) > 0$. Then:
    \[ K|_{\pi(\mathrm{Supp}(\mu) \cap \Lambda_0)} \equiv 0. \]
\end{lem}

\begin{proof}
Recall that $\Lambda_0 = \{\theta \in SM : \chi^+(\theta) = 0\}$ is invariant. Let $\tilde{\mu}$ denote the restriction of $\mu$ to $\Lambda_0 \cap \mathrm{Supp}(\mu)$. By the Birkhoff ergodic theorem, for $\tilde{\mu}$-a.e. $\theta$ there exists:
\begin{equation}\label{EQlem3.5}
    \Gamma_R^+(\theta) = \Gamma_R^-(\theta) =: \Gamma_R(\theta) = \widetilde{\mathrm{Ric}}(\theta) := \lim_{t \to +\infty} \frac{1}{t} \int_0^t \mathrm{Ric}(\phi^s(\theta))\,ds \in L^1(\tilde{\mu}) 
\end{equation}

with 
\[ \int_{\Lambda_0 \cap \mathrm{Supp}(\mu)} \widetilde{\mathrm{Ric}}(\theta)\,d\tilde{\mu} = \int_{\Lambda_0 \cap \mathrm{Supp}(\mu)} \mathrm{Ric}(\theta)\,d\mu. \]

From Lemma \ref{L1NEW}, and (\ref{EQlem3.5}) we have  $\widetilde{\mathrm{Ric}}(\theta) = 0$ holds $\tilde{\mu}$-a.e., implying:
\[ \int_{\Lambda_0 \cap \mathrm{Supp}(\mu)} \mathrm{Ric}(\theta)\,d\mu = 0. \]

\noindent Since $\chi^+(\theta) = \displaystyle\limsup_{t\to+\infty} \frac{1}{t} \int_0^t \mathrm{tr}(U_{\theta,u}(s))\,ds$ and $\mu(\Lambda_0) > 0$, we have $\widetilde{\mathrm{tr}\,U_u}(\theta) = 0$ $\mu$-a.e., where $\widetilde{\mathrm{tr}\,U_u}(\cdot)$ is the Birkhoff average of $\mathrm{tr}\,U_u(\cdot)$. 

\noindent Combining \eqref{EQ-Green} with the Birkhoff theorem and Cauchy-Schwarz yields:
$$
(n-1)\sqrt{-\int_{SM} \mathrm{Ric}(\theta)\,d\tilde{\mu}} = \int_{SM} \widetilde{\mathrm{tr}\,U_u}\,d\tilde{\mu} = \int_{SM} \mathrm{tr}\,U_u\,d\tilde{\mu} = 0.$$
\noindent Remark \ref{R_RIGI} now implies $R \equiv 0$ on $\mathrm{Supp}(\mu) \cap \Lambda_0$, as required.
\end{proof}

\noindent It is worth nothing that: In the case of dimension two and non-focal points (see (\cite{Barr-Pes}) for more details), the set $$\displaystyle\Lambda_0=\Big\{\theta: \limsup_{t\to\infty}\frac{1}{t}\int_0^tK(\gamma_{\theta}(s))ds=0\Big\},$$ 
where $K$ is the Gaussian curvature, captures directions where the average curvature along the geodesic vanishes asymptotically. In the case of non-positive curvature $$\Lambda_0=\{\theta: K(\gamma_{\theta}(t))=0, \,\, \text{for all}\,\, t\in \mathbb{R}\}.$$



\begin{proof}[\textbf{\emph{Proof of Theorem \ref{Theorem_Flat}}}]

It is clear that (c) implies both (a) and (b). We now prove the remaining implications:\\
\noindent (b) $\Rightarrow$ (c):\\
    Let $\mu$ be an invariant measure with $\operatorname{Supp} \mu = SM$ and $\mu(\Lambda_0) = 1$. By Lemma~\ref{L2'NEW}, we have
    \[
    0 \equiv K|_{\pi(\operatorname{Supp}(\mu))} \equiv K|_{\pi(SM)} \equiv K|_{M},
    \]
    which establishes (c).\\
    \noindent (a) $\Rightarrow$ (b):\\
    Under the given topological entropy condition, the variational principle and Pesin's formula imply that $\mathcal{L}(\Lambda_0) = 1$. Thus, the Liouville measure $\mathcal{L}$ satisfies (b).
\end{proof}

\section{Rigidity on Curvature Through Invariant Measures}
\begin{proof}[\textbf{\emph{Proof of Theorem \ref{Theorem_Rigidity}}}]
From \cite{Guimaraes}, since $M$ has curvature bounded from below and has no conjugate points, then the Ricci curvature is integrable for any invariant measure. 
    Assuming that $\mu$ is an invariant measure with $\mu(\Lambda_\alpha)=1$, then from (\ref{EQ-Form-Lyp-Exp})
$$\int_{SM}\text{tr}\,U_{u}\,d\mu=\int_{SM}\chi^{+}\,d\mu=\int_{\Lambda_{\alpha}}\chi^{+}\,d\mu =(n-1)\int_{\Lambda_{\alpha}}\alpha\, d\mu \geq (n-1)\sqrt{-\int_{SM}\emph{Ric}(\theta)\,d\mu}.$$

\noindent Note that for $\mu$-almost every point $\theta$, from  Birkhoff Ergodic Theorem, Remark \ref{RNEW1}, and (\ref{EQ-Form-Lyp-Exp}) both side of (\ref{EQ0NEW}) are converge to $\widetilde{\text{tr}\,U_u}(\theta)\leq (n-1)\sqrt{-\widetilde{\emph{Ric}}(\theta)}$, $\mu$ almost every $\theta$. Thus
\begin{eqnarray*}\label{EQ4NEW}
    \int_{SM}\text{tr}\,U_u \, d\mu &=& \int_{SM}\widetilde{\text{tr}\,U_u}(\theta)d\mu\leq (n-1)\int_{SM}\sqrt{-\widetilde{\emph{Ric}}(\theta)}d\mu\\
    &=&(n-1)\sqrt{-\int_{SM}\widetilde{\emph{Ric}}(\theta)\,d\mu}=(n-1)\sqrt{-\int_{SM}\emph{Ric}(\theta)\,d\mu}
\end{eqnarray*}
The two last inequalities provide 
\begin{equation*}\label{EQN1}
\int_{SM}\text{tr}\,U_{u}\,d\mu=(n-1)\sqrt{-\int_{SM}\emph{Ric}(\theta)\,d\mu}     
\end{equation*}

\noindent It follows from (\ref{EQ5NEW}) our result. 
\end{proof}

\begin{proof}[\textbf{\emph{Proof of Corollary \ref{C_Theorem_Rigidity}}}]
Note that if $\mu$ is ergodic with full support, then the invariance of $\Lambda_{\alpha}$ implies that $\mu(\Lambda_\alpha)=1$, then from Theorem \ref{Theorem_Rigidity}, $K|_{\pi(\text{Supp}(\mu)=SM}=\sqrt{-\int_{SM}\text{Ric}(\theta)d\mu}$.
\end{proof}

\begin{proof}[\textbf{\emph{Proof of Theorem \ref{Theorem_Rigidity2}}}]
   By proof of Theorem 1.4 of \cite{Nina2}, which relies on Kalinin's result (see \cite{Kal}), the ergodicity of the  Liouville measure $\mathcal{L}$ implies that all Lyapunov exponents are equal to $\alpha$. Consequently, $\mathcal{L}(\Lambda_\alpha)>0$ and the desired result follows directly from  Theorem \ref{Theorem_Rigidity}.
\end{proof}

\begin{proof}[\textbf{\emph{Proof of Theorem \ref{Theorem_Rigidity3}}}]
From (\ref{EQ-Form-Lyp-Exp}) we know that for every $\theta\in \text{Per}(\phi^t)$
$$\chi^{+}(\theta)=\frac{1}{\tau_{\theta}}\int_0^{\tau_\theta}\text{tr}\,U_{u}(\phi^s(\theta))ds=\int_{SM}\text{tr}\,U_{u}\,d\mu_{\theta},$$
where $\tau_\theta$ is the period of $\theta$. 
From hypothesis and Remark \ref{RNEW1}, the last inequality implies that 
$$\int_{SM}\text{tr}\,U_{u}\,d\mu_{\theta}=(n-1)\sqrt{-\int_{SM}\emph{Ric}(\theta)\,d\mu_{\theta}}.$$
Remark \ref{R_RIGI} provides that the sectional curvature $K_M$ is constant along $\pi(\phi^t(\theta))=\gamma_\theta(t)$, for every $\theta\in \text{Per}(\phi^t)$ and $U_u(\cdot)=\sqrt{-\int_{SM}\emph{Ric}(\theta)\,d\mu_{\theta}}\,I$.
The idea now is to use this property to extend the constant curvature to $\overline{\text{Per}(\phi^t)}$.\\
\ \\
\textbf{Claim:} Let $\theta=(x,v)\in \overline{\text{Per}(\phi^t)}$ and $w$ and unit vector in $T_xM$ linearly independent with $v$, then for every $t\in \mathbb{R}$
$$K_{x}(v,w)=K_{\gamma_{\theta}(t)}(\gamma'_{\theta}(t),W(t)),$$
where $W(t)$ is the parallel transport of $w$ along $\gamma_{\theta}(t)$.
\begin{proof}[\textbf{\emph{Proof of Claim}}]
There is a sequence of $\{\theta_n=(x_n,v_n)\}\subset  \text{Per}(\phi^t_{M})$  and a sequence $w_n\in T_{x_n}M$ such that $\theta_n$ converges to $\theta$ and $w_n$ converges to $w$. If $W_n(t)$ represents the parallel transport of $w_n$ along $\gamma_{\theta_n}(t)$, then since the sectional curvature is constant along the closed geodesic, then for every $t\in \mathbb{R}$ 
$$K_{x_n}(v_n,w_n)=K_{\gamma_{\theta_n}(t)}(\gamma'_{\theta_n}(t), W_n(t))$$
The left hand of the last equality converges to $K_{x}(v,w)$ and from the continuity of the parallel transport to the initial conditions, the right hand of the last equation converges to $K_{\gamma_{\theta}(t)}(\gamma'_{\theta}(t),W(t))$, as $n$ goes to infinite. Thus, we have the equality desired. 
\end{proof}

Let $\theta\in \overline{\text{Per}(\phi^t)}$, $\theta=(x,v)$ and $\{v,v_1, v_2, \dots,v_{n-1}\}$ an orthonormal basis of $T_xM$, and $\{\gamma'(t), v_1(t), v_2(t),\dots, v_{n-1}(t)\}$ a parallel orthonormal basic along $\gamma_{\theta}(t)$. So, from previous Claim 

\begin{eqnarray*}
    \text{Ric}(\phi^t_{M}(\theta))&=&\frac{1}{n-1}\sum_{i=1}^{n-1}R(\gamma'_{\theta}(t), v_i(t), \gamma'_{\theta}(t), v_i(t))\\
    &=&\frac{1}{n-1}\sum_{i=1}^{n-1} K_{\gamma_{\theta}(t)}(\gamma'_{\theta}(t),v_i(t))\\
    &=&\frac{1}{n-1}\sum_{i=1}^{n-1} K_x(v,v_i)=\text{Ric}(\theta).
\end{eqnarray*}
The last equality shows that the Ricci curvature $\text{Ric}|_{\overline{\text{Per}(\phi^t)}}$ is invariant for $\phi^t$, moreover, since it is a continuous function and $\phi^t|_{\overline{\text{Per}(\phi^t)}}$ is ergodic, then $\text{Ric}(\cdot)$ is a constant function, we said equal to $-\alpha^2$. Therefore, 
$U_u(\cdot)=\alpha\, I$ for the set of periodic orbits. This mean that $U_u(\cdot)=\alpha \, I$ in $\overline{\text{Per}(\phi^t)}$. The result follows from (\ref{Riccati}).
\end{proof}

\section{Rigidity of Regular Flow Equivalences}
\subsection{Rigidity for $1$-equivalences}

The following lemma serves a dual purpose: it is crucial for proving Theorems \ref{N-Conjugacy} and \ref{C1-Conjugacy}, while simultaneously generalizing \cite[Theorem 1.6]{Nina2} to a broader geometric setting.

\begin{lem}\label{Theorem_Not_Flat}
Let $M$ be a manifold with sectional curvature satisfying $-c^2 \leq K_M \leq -b^2 \leq 0$. If there exists a $\phi^t$-invariant measure $\mu$ with $\mu(\Lambda_\alpha) > 0$ for $\alpha \in \{b,c\}$, then:
\[ K_M|_{\pi(\mathrm{Supp}(\mu))} = -\alpha^2. \]
\end{lem}

\begin{proof}
The case $\alpha = c$ follows directly from Theorem \ref{Theorem_Rigidity}, since $\alpha = c \geq \sqrt{-\int_{SM} \mathrm{Ric}(\theta)d\mu}$. The case $\alpha = b = 0$ is handled by Lemma \ref{L2'NEW}. We thus assume $\alpha = b > 0$ with $\mu(\Lambda_b) = 1$.

For $\mu$-a.e. $\theta$, the Birkhoff ergodic theorem yields:
\[ b(n-1) = \lim_{t\to\infty} \frac{1}{t} \int_0^t \mathrm{tr}\,U_u(\phi^s(\theta))ds = \widetilde{\mathrm{tr}\,U_u}(\theta). \]

\noindent Note that $\displaystyle\int_{SM}\text{tr}\,U_{u}=\int_{SM}\widetilde{\text{tr}\,U_{u}}d\mu=b(n-1)$. Thus, since $K_M \leq -b^2$ implies $\mathrm{tr}\,U_u \geq (n-1)b$, we conclude $\mathrm{tr}\,U_u(\theta) = b(n-1)$ $\mu$-a.e. 
Finally, the Ricatti equation and symmetry of $U_{u}$ then force  
$\displaystyle \int_{SM}\text{Ric}(\theta)d\mu=-b^2$, which implies that $K_{M}|_{\pi(\text{Supp}(\mu))}=-b^2$. 
\end{proof}

\begin{proof}[\textbf{\emph{Proof of Theorem \ref{N-Conjugacy}}}]
Building on \cite{Guimaraes}, the curvature bound guarantees Ricci integrability. Let $h$ be a $1$-equivalence satisfying $h \circ \phi_M^t = \phi_N^{s(t)} \circ h$ with $s(t) \geq t$, $t\geq 0$.

For $\mu$-a.e. $\theta \in \Gamma_{M,N}$ (where $d_\theta h$ exists), there exist $C_1,C_2>0$ such that for all $t \in \mathbb{R}$ and $\eta \in G^u_\theta$:
\begin{equation}\label{EQ1}
C_1\|d_\theta\phi^t_M(\eta)\| \leq \|d_{h(\theta)}\phi^{s(t)}_N(d_\theta h(\eta))\| \leq C_2\|d_\theta\phi^t_M(\eta)\|.
\end{equation}

\noindent Note that, since $M$ has no conjugate points, then $\displaystyle\int_{SM}\text{Ric}\,d\mu\leq 0$.  Therefore,  If $\sup K_N\geq 0$, then $\displaystyle\int_{SM}\text{Ric}(\theta)\,d\mu\leq \sup K_N$. In this case, the equality only happens if both are equal to zero, however,  if $\displaystyle\int_{SM}\text{Ric}\,d\mu=0$, then $M$ is flat in $\pi(\text{Supp}(\mu))$ as we wished. (see \cite{Guimaraes}, \cite{Green}). Therefore, without loss of generality, we can assume that $\sup K_N<0$, which implies that $\phi^t_N$ is Anosov, since the curvature of $N$ is bounded from below,  therefore (\ref{EQ1}) implies 

\begin{equation}\label{EQ1*}
    \|d_{h(\theta)}\phi^{s(t)}_{N}(d_{\theta}h(\eta))\|\geq Ce^{\sqrt{-\sup K_N}s(t)}||d_{\theta}h(\eta)||.
\end{equation}
From (\ref{EQ1}) and (\ref{EQ1*}) 
\begin{equation}\label{EQ1**}
  \displaystyle \sqrt{-\sup K_N}\leq \limsup \frac{s(t)}{t}\sqrt{-\sup K_N}\leq \frac{\chi^{+}_{M}(\theta, \eta)}{n-1}. 
\end{equation}
The last inequality together with Remark \ref{RNEW1} and Birkhoff's ergodic theorem provides that, for $\mu$ almost every $\theta\in \Gamma_{M, N}$
 
\begin{equation}\label{EQ1***}
    \displaystyle (n-1)\sqrt{-\sup K_N}\leq \chi^{+}_{M}(\theta)\leq  (n-1)\sqrt{-\widetilde{\text{Ric}}\,(\theta)}.
\end{equation}

\noindent  Integrating with respect $\mu$ we have
\begin{eqnarray}\label{EQ2'New}
    (n-1)\sqrt{-\sup K_N}&\leq& \int_{SM}\chi^{+}(\theta)d\mu \leq (n-1)\int_{SM}\sqrt{-\tilde{Ric}\,(\theta)}d\mu \nonumber \\
    &\leq& (n-1)\sqrt{\int_{SM}-\widetilde{Ric}\,(\theta)d\mu}=(n-1)\sqrt{-\int_{SM}\text{Ric}(\theta)d\mu}.
\end{eqnarray}
Thus, we obtain that $$\sup K_N\geq \int_{SM}\text{Ric}(\theta)d\mu.$$
Note that, if we have equality in the last equation, then we have also equality in the equation (\ref{EQ2'New}), which implies that $\widetilde{\text{Ric}}\,(\theta)$ is constant $\mu$-almost everywhere, equal to $\displaystyle \int_{SM}\text{Ric}\,(\theta)d\mu$. Consequently, (\ref{EQ1***}) implies that for $\mu$-almost every $\theta$, $$\chi^{+}_{M}(\theta)=(n-1)\sqrt{-\int_{SM}\text{Ric}(\theta)d\mu}.$$
So, $\mu(\Lambda_{\alpha})=1$, for $\alpha= \sqrt{-\int_{SM}\text{Ric}(\theta)d\mu}$. Consequently, the Theorem \ref{Theorem_Rigidity} implies that $K_{M}|_{\pi(\text{Supp}(\mu))}=\displaystyle\int_{SM}\text{Ric}(\theta)d\mu=\sup K_N$.\\
To conclude the proof of the theorem we need to prove that $\displaystyle\lim_{t \to \infty}\frac{s(t)}{t}=1$. For this sake, note that the equation (\ref{EQ1**}) already imply that $\displaystyle\limsup_{t\to +\infty} \frac{s(t)}{t}=1$. Observe that, $\mu$-almost every $\theta$, from (\ref{EQ1}) and (\ref{EQ1*})
\begin{eqnarray}\label{AFTER-EQ1}
   \sqrt{-\int_{SM}\text{Ric}(\theta)d\mu}= \chi^{+}_{M}(\theta, \eta)&=&\displaystyle \lim_{t\to +\infty}\frac{1}{t}\log \|d_{\theta}\phi^{t}_{M}(\eta)\|=\displaystyle \liminf_{t\to +\infty}\frac{1}{t}\log \|d_{\theta}\phi^{t}_{M}(\eta)\| \nonumber\\
    &\geq & \liminf_{t\to +\infty}\frac{s(t)}{t}\sqrt{-\sup K_N}.
\end{eqnarray}
Thus $\displaystyle \liminf_{t\to +\infty}\frac{s(t)}{t}\leq 1$, but $s(t)\geq t$, then $\displaystyle\liminf_{t\to +\infty}\frac{s(t)}{t}=1$, as we wished.\\
Note that if $\displaystyle\liminf_{t\to +\infty}\frac{s(t)}{t}=1$, then (\ref{AFTER-EQ1}) implies that $\displaystyle\int_{SM}\text{Ric}(\theta)d\mu\leq \sup K_N $, and then $\displaystyle\int_{SM}\text{Ric}(\theta)d\mu = \sup K_N $, which completes the proof of the equivalences in the theorem. 
Assume now that $M$ and $N$ have the same dimension, our task is to prove that $K_N$ is also constant in $\pi(\text{Supp}(h_{*}\mu))$. For this sake, note that $K_M=-\alpha^2$, so (\ref{EQ1}) implies that $\chi^{+}_{N}(h(\theta))=(n-1)\alpha$, since $\sup K_N=-\alpha^2$, them using Lemma \ref{Theorem_Not_Flat} we conclude that $K_N|_{\pi(\text{Supp}(h_{*}(\mu)))}=-\alpha^2$.
\end{proof}
The techniques used in the proof of the previous theorem also work under weaker conditions. To understand these weak conditions, we introduce the following definition. 

\begin{defi}\label{E_1-measure}
    Two geodesic flows $\phi^{t}_{M}$ and $\phi^{t}_{N}$ have the property $\mathcal{E}_1$  along a measure $\mu$ if there is a map $h:U \to SN$ with $h \circ \phi^t_M = \phi^{s(t)}_{N} \circ h$ on $U$, and  $s(t)\geq t$ and satisfies \emph{(\ref{eq:bi_lipschitz})}.
\end{defi}
Using this definition on the hypotheses of Theorem \ref{N-Conjugacy}, we have the following result.
\begin{prop}\label{Prop1}
Theorem \ref{N-Conjugacy} remains valid when the global $\mathcal{E}_1$ condition is replaced by $\mathcal{E}_1$ along $\mu$.
\end{prop}

\begin{proof}[\textbf{\emph{Proof of Corolary \ref{N-Conjugacy-Liouville}}}]
It is enough apply the Theorem \ref{N-Conjugacy} for the Liouville measure $\mathcal{L}$, nothing that  $\mathcal{L}(SM\setminus \Gamma_{M,N})=1$ and $\text{Supp}(\mathcal{L})=SM$.
\end{proof}

\begin{proof}[\textbf{\emph{Proof of Theorem \ref{C1-Conjugacy}}}]
Since our equivalence is $C^1$, then $\Gamma_{N,M}=SM$. Therefore, for every $\theta\in\text{Per}(\phi^t_{M})$, $\mu_{\theta}(\Gamma_{M,N})=1$. So, similar arguments as in the proof of Theorem \ref{N-Conjugacy}, we obtain the condition $$\sup K_{N}|_{\mathcal{O}(h(\theta))}\geq  \int_{SM}\emph{Ric}(\theta)d\mathcal{\mu_{\theta}}.$$
Moreover, the equality holds if and only if $\displaystyle\lim_{t\to \infty}\frac{s(t)}{t}=1$, and 
\begin{equation}\label{EQN11}
K_{M}|_{\pi(\mathcal{O}(\theta))}= K_{M}|_{\pi(\text{Supp}(\mu_\theta))}=\int_{SM}\text{Ric}d\mu_{\theta}=\sup K_{N}|_{\mathcal{O}(h(\theta))}.  
\end{equation}
In other words, the section curvature is constant along periodic orbits. So, using the density of the periodic orbits and ergodicity of geodesic flow, a similar argument  to the proof Theorem \ref{Theorem_Rigidity3} allows us to conclude that $\text{Ric}(\cdot)$ is a constant function, we said equal to $-\alpha^2$. Therefore, from (\ref{EQN11}) we obtain that 
$$K_{M}|_{\pi(\mathcal{O}(\theta))}=-\alpha^2, \,\,\, \text{for all}\,\,\, \theta\in \text{Per}(\phi^t_{M}).$$
Consequently, from \cite{Nina} we have density of $\text{Per}(\phi^t_{M})$, therefore, the continuity of $K_M$ allows us to conclude that $K_M=-\alpha^2$ as we wished.\\
When $M$ and $N$ have the same dimension the proof follows the same lines like the end of the proof of Theorem \ref{N-Conjugacy}, from (\ref{EQ1}), $\chi^{+}_{N}(h(\theta))=(n-1)\alpha$ and $K_M=-\alpha^2$.
\end{proof}

\begin{remark}
 In the proof of the previous theorem, the Anosov condition of $\phi^t_{M}$ can be replaced by $M$ without conjugate points, $\phi^t_{M}$ ergodic, and density of periodic orbits. This condition is satisfied for a wide range of manifolds.
\end{remark}
The conditions on density of periodic orbits and ergodicity are useful to spread the constant curvature throughout the whole manifold. However, the following proposition can be obtain from the proof for Theorem \ref{C1-Conjugacy} and holds for each geodesic. 
\begin{prop}\label{Prop_2} 
 Let $M$ and $N$ be two complete Riemannian manifolds whose geodesic flows $\varphi_{M}^{t}$ and $\varphi_{N}^{t}$ satisfies the property $\mathcal{E}_1$  with $C^1$-equivalence $h$ along a periodic measure $\mu_{\theta}$. Assume that $M$ has no conjugate points along $\gamma_{\theta}(t)$, then  $$ K_{N}|_{\mathcal{O}(h(\theta))}\geq  \int_{SM}\emph{Ric}(\theta)d\mathcal{\mu_{\theta}}.$$ 
 The equality holds  if and only if, $M$ has constant negative sectional curvature along $\gamma_{\theta}(t)$ and $\displaystyle\lim_{t\to +\infty}\frac{s(t)}{t}=1$. Moreover, if $M$ and $N$ have the same dimension, then $N$ has constant negative sectional in $\gamma_{h(\theta)}(t)$.
\end{prop}

\subsection{More General Rigidity}
This section extends Theorems \ref{N-Conjugacy} and \ref{C1-Conjugacy} to equivalences where the reparametrization $s(t)$ satisfies $s(t) \geq at$ for some $a > 0$, using conformal metric deformations.
\subsubsection{No conjugate points and homotheties}
We investigate relationships between geodesic flows on the same manifold under different metrics, particularly homothetic transformations.
Let $(M,g)$ be a Riemannian manifold and consider the conformal metrics $g_r = e^{2r}g$ for $r \in \mathbb{R}$. Denote by $M_r$ the manifold $M$ with metric $g_r$. Key observations (see \cite{doCarmo} for more details):

\begin{itemize}
    \item The Levi-Civita connections satisfy $\nabla_r = \nabla$.
    \item Sectional curvatures relate by $K_r(v,w) = e^{-2r}K(v,w)$ for linearly independent $v,w \in T_xM$.
    \item Geodesics coincide as sets, but velocities scale: $\|\gamma'(t)\|_r = e^r\|\gamma'(t)\|$.
    \item Jacobi fields satisfy $\|J(t)\|_r = e^r\|J(t)\|$ .
\end{itemize}

\begin{lem}\label{Homo-NCP}
$(M,g_r)$ has no conjugate points if and only if $(M,g)$ has no conjugate points.
\end{lem}
\begin{proof}
Immediate from the norm scaling of Jacobi fields.
\end{proof}

\begin{lem}\label{lemma_Homothe}
For any complete $(M,g)$, there exists an orbit-isometric $C^1$-conjugacy $h:SM \to SM_r$ between $\phi_M^t$ and $\phi_{M_r}^t$.
\end{lem}
\begin{proof}
Define $h(x,v) = (x,e^{-r}v)$ and $s(t) = e^rt$. For $\theta = (x,v) \in SM$ and $\widetilde{\theta} = h(\theta)$, the geodesics satisfy $\gamma_{\widetilde{\theta}}(t) = \gamma_\theta(e^{-r}t)$. By definition of $h$ and the geodesic flows $\varphi^{t}_{M}$ and $\varphi^{t}_{M_{r}}$ we have that for all $t\in \mathbb{R}$
\begin{eqnarray*}
h\circ\phi^{s^{-1}(t)}_{M}(\theta)&=&h\circ\phi^{e^{-r}t}_{M}(\theta)
=h(\gamma_{\theta}(e^{-r}t),\gamma'_{\theta}(e^{-r}t))\\
&=&(\gamma_{\theta}(e^{-r}t),e^{-r}\gamma'_{\theta}(e^{-r}t))
=(\gamma_{\widetilde{\theta}}(t),\gamma'_{\widetilde{\theta}}(t))\\
&=&\phi^{t}_{M_{r}}(x,e^{-r}t)
=\phi^{t}_{M_{r}}\circ h(\theta).
\end{eqnarray*}

\noindent Hence $h\circ\phi^{t}_{M}=\phi^{s(t)}_{M_{r}}\circ h$, for all $t\in\mathbb{R}$.
\end{proof}

Observe that, if $r<0$, then the last lemma implies that $\phi^{t}_{M}$ and $\phi^{t}_{M_{r}}$ have the property $\mathcal{E}_1$.
\subsubsection{Generalized Equivalence Properties}
\begin{defi}
    Geodesic flows $\phi^{t}_{M}$ and $\phi^{t}_{N}$ have the property $\mathcal{E}_a$ $(a>0)$, if they are $1$-equivalent  reparametrization $s(t)$ satisfying $s(t)\geq at$, $t\geq 0$.
\end{defi}

 \begin{teo}\label{N-Conjugacy-a} Let $M$ and $N$ be two complete Riemannian manifolds whose geodesic flows  $\phi_{M}^{t}$ and $\phi_{N}^{t}$ satisfies the property $\mathcal{E}_a$ with $1$-equivalence $h$. If $M$ has no conjugate points and there exists an invariant measure $\mu$ with $\mu(\Gamma_{M,N}) = 1$, then: $$a^2\sup K_{N}\geq\int_{SM}\emph{Ric}(\theta)d\mu.$$ The equality holds if and only if  $K_{M}|_{\pi(\text{Supp}(\mu))}=a^2\sup K_N$ and $\displaystyle\lim_{t\to +\infty}\frac{s(t)}{t}=a$. \\
 Moreover, when $\dim M = \dim N$, $N$ has constant negative curvature on $\pi(\mathrm{Supp}(h_*\mu))$.
\end{teo}
\begin{proof}
Since $a>0$ there is $r_0$ such that $a=e^{-r_0}$. Consider the manifold $(N_{r_0}, g_{r_0})$, then 
$$a^2K_N= K_{N_{r_0}}.$$
From Lemma \ref{lemma_Homothe} there is an isometry $h$  such that $h\circ\phi^{t}_{N}=\phi^{u(t)}_{N_{r_0}}\circ h$, where $u(t)=e^{r_0}t$.\\
Also, by hypotheses, there is an $1$-equivalence $f$ with $f \circ \phi^t_{M}=\phi ^{s(t)}_{N}\circ f$, where $s(t)\geq at$, $t\geq 0$. 

$$h\circ f \circ \phi^t_{M} = h \circ \phi^{s(t)}_{N}\circ f= \phi_{N_{r_0}}^{u(s(t))}\circ h\circ f.$$
Note that $u(s(t))=e^{r_0}s(t)\geq e^{r_0}at=t$. Moreover, since $h$ is an isometry, then $h\circ f$ is a $1$-equivalence. Consequently, we obtain our result from Theorem \ref{N-Conjugacy}.
\end{proof}
A similar result to Theorem \ref{Theorem_Rigidity3} is also obtained with the property $\mathcal{E}_{a}$. Its proof follows the same lines as the proof of Theorem \ref{N-Conjugacy-a} but using Theorem \ref{C1-Conjugacy} instead of Theorem \ref{N-Conjugacy}.

 \begin{teo}\label{Theorem_C^1-Conjugacy-a}

 Let $M$ and $N$ be two complete Riemannian manifolds whose geodesic flows $\phi_{M}^{t}$ and $\phi_{N}^{t}$ satisfy the property $\mathcal{E}_a$ with $C^1$-equivalence $h$. If $M$ has finite volume and $\phi^t_{M}$ is Anosov, then:  $$ a^2\,K_{N}|_{\mathcal{O}(h(\theta))}\geq  \int_{SM}\emph{Ric}(\theta)d\mathcal{\mu_{\theta}}, \, \, \, \text{for all} \, \, \, \theta \in \text{Per}(\phi^t_{M}).$$ 
  The equality holds for all geodesic if and only if, $M$ has constant negative sectional curvature and $\displaystyle\lim_{t\to +\infty}\frac{s(t)}{t}=a$. 
  Moreover, when $\dim M = \dim N$, $N$ has constant negative curvature on $\pi(\mathrm{Supp}(h_*\mu))$.
\end{teo}

\begin{remark}
The property $\mathcal{E}_a$ admits a natural measure-theoretic extension $\mathcal{E}_a(\mu)$ for invariant measures $\mu$, producing corresponding generalizations of \emph{Theorems \ref{N-Conjugacy-a} and \ref{Theorem_C^1-Conjugacy-a}} in complete analogy with \emph{Propositions \ref{Prop1} and \ref{Prop_2}}.
\end{remark}

\section{Some Applications}\label{Some Applications}

\subsection{Generalization of Some Rigidity Results}
In this section we focus on made a simple proof of some rigidity result at \cite[Theorem 1.1, Theorem 1.2 , and Theorem 1.4]{IR2}.
\begin{teo}[\cite{IR2}, Theorem 1.1]
 If $M$ is a complete manifold of finite volume with sectional curvature bounded from below by $-c^2$ and Anosov geodesic flow. If the constant of contraction of the flow $\lambda=e^{-c}$, then $K_M=-c^2$.
\end{teo}
 \begin{proof} 
The curvature condition implies (see \cite[Lemma 2.16]{Kn}) the following uniform growth estimates for unstable vectors:
\begin{equation*}\label{eq:growth_estimate}
C\lambda^{-t}\|\eta\| \leq \|d_{\theta}\varphi^{t}_{M}(\eta)\| \leq e^{ct}\|\eta\|\sqrt{1+c^2}, \quad t \geq 0, \,\, \eta \in E^u(\theta).
\end{equation*}
Taking $\lambda = e^{-c}$, we obtain $\chi^+(\theta) = (n-1)c$ for all $\theta \in SM$, which implies $\mathcal{L}(\Lambda_c) = 1$ for the Liouville measure $\mathcal{L}$. \\
Since the geodesic flow is Anosov and $M$ has finite volume, the Liouville measure is ergodic. Furthermore, the inequality 
\[
c \geq \sqrt{-\int_{SM} \text{Ric}(\theta)\,d\mathcal{L}}
\]
holds, and thus the conclusion follows directly from Corollary \ref{C_Theorem_Rigidity}.
 \end{proof}

Our next result extends \cite[Theorem 1.2 and Theorem 1.4]{IR2} in two significant aspects: the manifolds need not be compact, nor is an Anosov geodesic flow required.

\begin{teo}\label{thm:curvature_comparison}
Let $M$ and $N$ be Riemannian manifolds whose geodesic flows are $1$-equivalent with $s(t) = t$. Assume that:
\begin{enumerate}
    \item[\emph{1}.] $M$ has finite volume;
    \item[\emph{2}.]  $M$ has curvature bounded below and no conjugate points;
    \item[\emph{3}.]  $\inf K_{M} \geq \sup K_N$.
\end{enumerate}
Then $K_M \equiv \sup K_N$ on $M$. Moreover, if $\dim M = \dim N$, then $K_M = K_N \equiv \sup K_N$.
\end{teo}

\begin{proof}
Since $\phi_M^t$ and $\phi_N^t$ satisfy property $\mathcal{E}_1$, Theorem \ref{N-Conjugacy} yields
\[ \sup K_{N} \geq \int_{SM} \text{Ric}(\theta)\, d\mathcal{L}. \]

The hypothesis $\inf K_{M} \geq \sup K_N$ implies that $\text{Ric}(\theta) \equiv \inf K_M$.
Consequently,
\[ \sup K_{N} = \int_{SM} \text{Ric}(\theta)\, d\mathcal{L}. \]

\noindent The conclusion follows directly from Corollary \ref{N-Conjugacy-Liouville}.
\end{proof}

\textbf{Acknowledgements:} 
The author is deeply grateful to Nestor Nina Zarate and Alexander Cantoral for their invaluable insights and stimulating mathematical discussions during the preparation of this work. Their contributions significantly enriched this research. Special thanks are extended to all members of the School of Mathematics (Zhuhai) at Sun Yan-sen University for providing such an exceptional academic environment that greatly supported this research.

\bibliographystyle{plain}
\pagestyle{empty}


\noindent Sergio Roma\~{n}a\\
School of Mathematics (Zhuhai), Sun Yat-sen University, 519802, China\\
sergio@mail.sysu.edu.cn

\end{document}